\documentclass[11pt,a4paper]{article}

\usepackage[cp1251]{inputenc}
\usepackage[english,russian]{babel}

\usepackage{textcomp}
\usepackage{amssymb,amsmath,amscd,amsthm}
\usepackage{array}
\usepackage{cite}

\makeatletter  %%%%%%%%%%%%%%%%  включение работы с командами содержащими символ @

\def\@sect#1#2#3#4#5#6[#7]#8{\ifnum #2>\c@secnumdepth
     \let\@svsec\@empty\else
     \refstepcounter{#1}\edef\@svsec{\csname the#1\endcsname.\hskip 0.4em}\fi
     \@tempskipa #5\relax
      \ifdim \@tempskipa>\z@
        \begingroup #6\relax
          \@hangfrom{\hskip #3\relax\@svsec}{\interlinepenalty \@M #8\par}%
        \endgroup
       \csname #1mark\endcsname{#7}\addcontentsline
         {toc}{#1}{\ifnum #2>\c@secnumdepth \else
                      \protect\numberline{\csname the#1\endcsname.}\fi
                    #7}\else
        \def\@svsechd{#6\hskip #3\relax  %% \relax added 2 May 90
                   \@svsec #8\csname #1mark\endcsname
                      {#7}\addcontentsline
                           {toc}{#1}{\ifnum #2>\c@secnumdepth \else
                             \protect\numberline{\csname the#1\endcsname.}\fi
                       #7}}\fi
     \@xsect{#5}\noindent\vspace{-2ex}}

\@addtoreset {equation}{section}
\makeatother  %%%%%%%%%%%%%%%%  отключение работы с командами содержащими символ @

\newcommand{\beq}{\begin{equation}}
\newcommand{\beqnt}{\begin{equation}\notag}
\newcommand{\eeq}{\end{equation}}
\newcommand{\bml}{\begin{multline}}
\newcommand{\eml}{\end{multline}}
\newcommand{\bmlnt}{\begin{multline*}}
\newcommand{\emlnt}{\end{multline*}}
\newcommand{\bpr}{\begin{proof}}
\newcommand{\bpra}[1]{\begin{proof}[#1]}
\newcommand{\epr}{\end{proof}}
\newcommand{\fref}[1]{{\rm(\ref{#1})}}                      %%%%%%%%%% formula ref      %%%%%
\theoremstyle{definition}
\newtheorem{teoRU}{Теорема}[section]
\newtheorem{lemRU}{Лемма}[section]

\newtheorem{corRU}{Следствие}[section]

\newtheorem{remRU}{Замечание}[section]

\newtheorem{exaRU}{П\ р\ и\ м\ е\ р\ }[section]

\newcommand{\eissue}{\vfill{\tiny Compiled {\tt \number\day.\number\month.\number\year} from {\tt {\jobname}.tex}}}
\newcommand{\rissue}{\vfill{\tiny Получено \today\ из файла {\tt {\jobname}.tex}}}

\newcommand{\UB}{{\mathbf U}}%%%%%%% всякие множества управлений
\newcommand{\VB}{{\mathbf V}}%%%%%%% всякие множества помехи

\newcommand{\RA}{{\mathbb R}}%%%%%%%%%  действительные числа
\newcommand{\NA}{{\mathbb N}}%%%%%%%%%  натуральные числа

\newcommand{\barf}[1]{\hat #1}%%%%%   обозначение <расширенния> многозначного отображения 

\newcommand{\Fg}[2]{{#1}_{#2}}%%%%%   многозначное отображение отвечающее игре двух и более лиц

\newcommand\argmin{\operatornamewithlimits{\mathrm{argmin}}}
\newcommand\argmax{\operatornamewithlimits{\mathrm{argmax}}}

\newcommand\diametr{\operatorname{\mathbf{diam}}}%%%%%%%%%%%%%%%% диаметр множества

\newcommand\di{\operatorname{\rho}}  %%%%%%%%%%%%%%%%%%%%%  метрика
\newcommand\dis{\operatorname{\mathbf d}}  %%%%%%%%%%%%%%%%%%%%%  метрика

\newcommand{\nint}[2]{{{#1}..{#2}}}
\newcommand{\mydef}{\operatorname{:\!=}}%%%%%%%%%%%%%%%%%%%%% знак <<равно по определению>> конструкции
\newcommand{\myinf}{\operatornamewithlimits{\inf\vphantom{\sup}}}%%%%%%%%%%%%%%%%%%%%% знак inf выравненый <<под sup>>
\newcommand{\fix}[1]{{\mathbf{Fix}({#1})}}%%%%%   множество неподвижных точек отображения
\newcommand{\Ox}{O(X)}%%%%%   конечные покрытия X
\newcommand{\Ofo}{O_\mathrm{fo}(X)}%%%%%   конечные открытые покрытия X
\newcommand{\Ofc}{O_\mathrm{fc}(X)}%%%%%   конечные замкнутые покрытия X
\newcommand{\Ofcd}[1]{O^{#1}_\mathrm{fc}(X)}%%%%%   конечные замкнутые покрытия X диаметра \delta
\newcommand{\Ofod}[1]{O^{#1}_\mathrm{fo}(X)}%%%%%   конечные открытые покрытия X диаметра \delta
\newcommand{\poc}{\mathrel\sqsubseteq}%%%%%  отношение вписанности на покрытияx

\renewcommand{\ge}{\geqslant}                       %%%%%%%%%%%%%%% правильные знаки "Больше равно", "Меньше равно"
\renewcommand{\le}{\leqslant}                       %%%%%%%%%%%%%%% правильные знаки "Больше равно", "Меньше равно"
\newcommand{\myle}{\mathrel\preccurlyeq}                %%%%%%%%%%%%%%% знак частичного порядка 
\newcommand{\myand}{\ensuremath{\land}}   %%%%%%%%%%%%%%%%%%%%% знак <<и>>
\newcommand{\myimp}{\ensuremath{\Rightarrow}}%%%%%%%%%%%%%%%%%%%%% знак импликации
\newcommand{\myeqv}{\ensuremath{\Leftrightarrow}}%%%%%%%%%%%%%%%%%%%%% знак эквиваленции
\newcommand{\myll}{\ensuremath{\forall}}         %%%%%%%%%%%%%%%%%%%%% квантор всеобщности
\newcommand{\myxst}{\ensuremath{\exists}}   %%%%%%%%%%%%%%%%%%%%% квантор существования
\title{О неподвижных точках многозначных отображений}  % Declares the document's title.

\author{Д.\,А.\,Серков}      % Declares the author's name.

\date{Институт математики и механики им. Н.Н.Красовского УрО РАН\\Уральский федеральный университет им. Б.Н.Ельцина\\ Екатеринбург, Россия\\e-mail: serkov@imm.uran.ru}

\usepackage[unicode,pdfpagemode=None,linkcolor=black,citecolor=black,colorlinks]{hyperref}%%%%,pagebackref

\begin{document}

\maketitle

%\tableofcontents

\abstract{В работе рассмотрены условия существования неподвижных точек многозначных отображений, не опирающиеся на линейную структуру множества. Даны описания множеств неподвижных точек для отображений с замкнутым графиком --- в компактных Хаусдорфовых пространствах, в метризуемых пространствах, а также для непрерывных функций --- в пространствах со сходимостью (топология Фреше) и с топологией Скотта. Даны приложения к задачам о равновесиях --- приведены описания седловых точек и равновесий по Нэшу для компактных Хаусдорфовых и метризуемых пространств чистых стратегий игроков.

The paper discusses the conditions for the existence of fixed points of multivalued mappings that are not based on the linear structure of the set. The descriptions of the sets of fixed points for mappings with closed graph in compact Hausdorff spaces and in metrizable spaces, as well as for continuous functions in spaces with convergence (Frechet topology) and in spaces with Scott topology are provided. Applications to the problem of the equilibrium are given: the sets of saddle points and of Nash equilibria for compact Hausdorff and metrizable spaces of strategies of players are described. 
}

%%%%%%%%%%%%%%%%%%%%%%%%%%%%%%%%%%%%%%%%%%%%%%%%%%%%%%%%%%%%%%
\section*{Введение}
\addcontentsline{toc}{section}{Введение}

Рассмотрим автоморфизм $F$ множества $X$ как динамику некоторой системы с дискретным временем и фазовым пространством $X$. В случае компактности $X$ всякая траектория такой системы имеет предельные точки. Эти предельные точки либо образуют цикл, либо являются стационарными точками (циклами длины 1) и, как следствие,   неподвижными точками $F$. 

Многие достаточные условия существования неподвижных точек можно рассматривать как условия, препятствующие образованию циклов, состоящих из нескольких точек:

--- <<принцип сжатия>>: расстояние между образами \emph{произвольных} двух точек меньше, чем расстояние между этими точками (цикл обязан стягиваться в точку);

--- <<принцип направленности>>: множество $X$ частично упорядоченно, а  отображение $F$ изотонно по отношению к этому порядку (возможны лишь одноточечные циклы);

Идея предлагаемых условий состоит в ограничении размеров циклов: цикл из одной точки обязан содержаться в любом множестве, содержащем эту точку.

%%%%%%%%%%%%%%%%%%%%%%%%%%%%%%%%%%%%%%%%%%%%%%%%%%%%%%
\section{Определения}

1. Для произвольного множества $X$ обозначим $2^X$ множество всех подмножеств $X$. Пусть имеется не пустое множество $X$ и многозначное отображение $X\ni x\mapsto F(x)\in2^{X}$. 
Обозначим
$$
\fix{F}\mydef\{x\in X\mid x\in F(x)\}
$$
--- множество всех неподвижных точек отображения $F$. 

По заданному отбражениею $F$ определим отображение $2^X\ni Y\mapsto\barf{F}(Y)\in2^X$ следующим образом:
\beqnt
\barf{F}(Y)\mydef\Biggl(\bigcup_{y\in Y}F(y)\Biggr)\bigcap Y=\bigcup_{y\in Y}F(y)\cap Y.
\eeq
Заметим, что множество $\barf{F}(Y)$ является множеством <<кандидатов>> из множества $Y$ во множество неподвижных точек: легко проверить <<от противного>>, что  элементы из $Y$, не попавшие в $\barf{F}(Y)$, заведомо не принадлежат $\fix{F}$. 

В частности, для любого синглетона $S_x\in2^X$, $x\in X$: 
$$
S_x\mydef\{x\},
$$
множество  $\barf{F}(S_x)$ не пусто, тогда и только тогда, когда $x\in\fix{F}$.

Непосредственно из определения получим: для всяких $Y,Y'\in2^X$ выполняются соотношения
\beq\label{Lf1-inc1}
\barf{F}(Y)\subset Y,
\eeq
\beq\label{Lf1-inc2}
Y'\subset Y\quad\Rightarrow\quad\barf{F}(Y')\subset \barf{F}(Y),
\eeq
\beq\label{Lf1-inc3}
\fix{F}\cap Y=\fix{F}\cap\barf{F}(Y).
\eeq

Для произвольного семейста  $(Z_\tau)_T\subset2^X$ из определения $\barf{F}$ и приведенных соотношений получаем
\beq\label{Lf2-inc1}
 \bigcup_{\tau\in T}\barf{F}(Z_\tau)\subset \bigcup_{\tau\in T}Z_\tau,
\eeq
\beq\label{Lf2-inc2}
\fix{F}\bigcap \bigcup_{\tau\in T}Z_\tau=\fix{F}\bigcap \bigcup_{\tau\in T}\barf{F}(Z_\tau),
\eeq

Здесь включение \fref{Lf2-inc1} следует из включения \fref{Lf1-inc1}, включение \fref{Lf2-inc2} --- из \fref{Lf1-inc3}.

2. На подмножествах из $2^X$ определим отношение \emph{быть вписанным}, которое будем обозначать $\poc$: для произвольных
$(Z_\tau)_{T},(Z'_{\tau'})_{T'}\subset2^X$ будем говорить, что семейство $(Z'_{\tau'})_{T'}$ вписано в семейство $(Z_\tau)_{T}$ и обозначать это $(Z'_{\tau'})_{T'}\poc (Z_\tau)_{T}$, если для произвольного $\tau'\in T'$ найдется $\tau\in T$ такой, что $Z'_{\tau'}\subset Z_{\tau}$. 

Такая структура порядка в дальнейшем будет полезна ввиду следующего свойства отображения $\barf{F}$:
из импликации \fref{Lf1-inc2} следует, что для любых $(Z'_{\tau'})_{T'},(Z_{\tau})_{T}\subset2^X$, если $(Z'_{\tau'})_{T'}\poc (Z_\tau)_{T}$, то выполняется включение
\beq\label{Lf2-inc4}
\bigcup_{\tau'\in T'}\barf{F}(Z'_{\tau'})\subset \bigcup_{\tau\in T}\barf{F}(Z_\tau).
\eeq

Обозначим $\Ox$ --- множество всех покрытий $X$, то есть всех подмножеств $(O_\iota)_I\subset2^X$ таких, что $\cup_{\iota\in I}O_\iota=X$. 

Пусть $\tau(X)$ --- некоторая топология в $X$ (множество всех открытых множеств). 
Обозначим  $\Ofo$ $(\Ofc)$--- множество всех конечных открытых (замкнутых) покрытий $X$.

%%%%%%%%%%%%%%%%%%%%%%%%%%%%%%%%%%%%%%%%%%%%%%%%%%%%%%
\section{Критерий существования неподвижных точек}

Множество неподвижных точек отображения $F$ можно характеризовать следующим образом.

\begin{teoRU}\label{LOx} Для любого множества $X\neq\varnothing$ и любого отображения $F:X\mapsto2^X$ справедливо равенство
\beq\label{fix-eq-nocycle-x}
\fix{F}=\bigcap_{(O_\iota)_{ I}\in \Ox} \bigcup_{\iota\in I}\barf{F}(O_\iota).
\eeq
\end{teoRU}

\begin{remRU}
Так как среди покрытий имеется наименьший элемент по отношению вписанности --- покрытие состоящее из синглетонов множества $X$:
$(S_x)_X$, то для этого покрытия и любого $(O_\iota)_I\in \Ox$ в силу \fref{Lf2-inc4} выполнены включения 
$$
\bigcup_{x\in X}\barf{F}(S_x)\subset\bigcup_{\iota\in I}\barf{F}(O_\iota),
$$
а значит, и равенство
\beq\label{fix-eq-singleton}
\fix{F}=\bigcup_{x\in X}\barf{F}(S_x).
\eeq

Однако, нам важн\'о именно соотношение \fref{fix-eq-nocycle-x}, так как оно открывает возможности для использования свойств множества $X$ и отображения $F$.
\end{remRU}

В основе приводимых ниже следствий  равенства \fref{fix-eq-nocycle-x} лежит тот интуитивно очевидный факт, что если некоторое замкнутое подмножество $G$ компактного множества $X\times X$ имеет точки сколь угодно близкие к диагонали $\{(x,x)\mid x\in X\}$, то $G$ содержит и некоторый элемент $(\bar x,\bar x)$ диагонали. Если же в качестве множества $G$ рассматривать график многозначного отображения $F$, то элемент $\bar x$ оказывается неподвижной точкой $F$. 

\begin{teoRU}\label{TfO}
Пусть $X$ --- компактное хаусдорфово пространство, а отображение $F$ имеет зам\-кнутый график. Тогда 
\beq\label{fix-eq-nocycle-o}
\fix{F}=\bigcap_{(O_\iota)_{ I}\in \Ofo} \bigcup_{\iota\in I}\barf{F}(O_\iota)=\bigcap_{(O_\iota)_{ I}\in \Ofc} \bigcup_{\iota\in I}\barf{F}(O_\iota).
\eeq
В частности, множество $\fix{F}$ не пусто, если и только если выполнено условие
\beq\label{con-cov-not-empt}
(\forall(O_\iota)_{ I}\in \Ofc)(\exists\bar\iota\in I)\qquad\barf{F}(O_{\bar\iota})\neq\varnothing.
\eeq
\end{teoRU}

\begin{corRU}\label{TfTS}
Пусть $X$ --- компактное хаусдорфово пространство, а функция $f:X\mapsto X$ непрерывна. Тогда 
\beq\label{fix-eq-Ocf-singl}
\fix{f}=\bigcap_{(O_\iota)_{ I}\in \Ofo} \bigcup_{\iota\in I}\barf{f}(O_\iota)=\bigcap_{(O_\iota)_{ I}\in \Ofc} \bigcup_{\iota\in I}\barf{f}(O_\iota).
\eeq
В частности, функция $f$ имеет неподвижную точку тогда и только тогда, когда 
\beq\label{con-cov-not-empt-f}
(\forall(O_\iota)_{ I}\in \Ofc)(\exists\bar\iota\in I)\qquad\barf{f}(O_{\bar\iota})\neq\varnothing.
\eeq
\end{corRU}

Пусть $X$ --- метрическое пространство с метрикой $\di:X^2\mapsto[0,\infty)$. 
Расстояние от точки $x\in X$ до множества $A\subset X$, порожденное метрикой $\di$, будем обозначать $\dis(x,A)$:
$$
\dis(x,A)\mydef\inf_{a\in A}\di(a,x).
$$

Для всякого $Y\in2^X$ обозначим $\diametr(Y)\mydef\sup_{y,y'\in Y}\di(y,y')$ --- диаметр множества $Y$. 
Для всякого $\delta>0$ через $\Ofcd{\delta}$ ($\Ofod{\delta}$) обозначим подмножество  $\Ofc$ замкнутых ($\Ofo$ открытых) конечных порытий $(O_\iota)_I$ диаметра не более $\delta$: 
$$
\max_{\iota\in I}\diametr(O_\iota)\le\delta.
$$
Напомним, что для компактного мерического пространства $X$ при любом $\delta>0$ множества $\Ofcd{\delta}$, $\Ofod{\delta}$ не пусты.

\begin{teoRU}\label{Tf2}
Пусть $X$ --- метрический компакт, а отображение $F$ имеет замкнутый график. Тогда для произвольного семейства покрытий
\beq\label{dim-diam-cov}
(O_{k\iota})_{I_k}\in \Ofcd{\delta_k}\cup\Ofod{\delta_k},\quad\delta_k>0, \quad k\in\NA,\qquad \lim_{k\to\infty}\delta_k=0
\eeq
выполнено равенство
\beq\label{fix-eq-nocycle-2}
\fix{F}=\bigcap_{k\in\NA} \bigcup_{\iota\in I_k}\barf{F}(O_{k\iota}).
\eeq
В частности, множество $\fix{F}$ не пусто тогда и только тогда, когда 
\beq\label{x-close-to-Fx}
(\forall \delta>0)(\exists x_\delta\in X)\qquad  \dis(x_\delta,F(x_\delta))\le\delta.
\eeq
\end{teoRU}

\begin{corRU}\label{Tf2S}
Пусть $X$ --- метрический компакт, а функция $f$ непрерывна. Тогда для произвольного семейства покрытий вида \fref{dim-diam-cov}
выполнено равенство
\beq\label{fix-eq-nocycle-2-S}
\fix{f}=\bigcap_{k\in\NA} \bigcup_{\iota\in I_k}\barf{f}(O_{k\iota}).
\eeq
В частности, функция $f$ имеет неподвижную точку тогда и только тогда, когда 
\beq\label{x-close-to-fx}
(\forall \delta>0)(\exists x_\delta\in X)\qquad  \di(x_\delta,f(x_\delta))\le\delta.
\eeq
\end{corRU}

%%%%%%%%%%%%%%%%%%%%%%%%%%%%%%%%%%%%%%%%%%%%%%%%%%%%%%
\subsection{Пример}\label{exaTf2}

Приведем по возможности простой пример использования теоремы \ref{Tf2}. В качестве метрического компакта $X$ возьмем отрезок $[0,1]$ с естественной метрикой: $\di(x,y)\mydef|x-y|$. Выберем и зафиксируем бесконечную дизъюнктную последовательность $([a_i,b_i])_\NA$ отрезков из $X$ ненулевой длины. Положим  
$$
G(x)\mydef\bigcup_{i\in\NA} G_{[a_i,b_i]}(x),
\qquad
G_{[a,b]}(x)\mydef
\begin{cases}
\argmax_{y\in[a,b]}|y-x|,&x\in[a,b],\\
\varnothing,&x\notin[a,b].
\end{cases}
$$
Определим отображение $F:X\mapsto2^{X}$ следующим образом: графиком отображения $F$ является замыкание в $\RA^2$ графика отображения $G$. Заметим, что при некоторых $x\in X$ значения отображения $F$ являются не выпуклыми или пустыми.

Легко видеть, что отображение $F$ удовлетворяет условиям теоремы \ref{Tf2}: в силу требования дизъюнкности длины отрезков $[a_i,b_i]$ по необходимости сремятся к нулю при $i\to\infty$ и, значит,  для любого $\delta>0$ при достаточно большом значении $i$ будет выполняться неравенство $\dis(a_i,F(a_i))\le\delta$. Cледовательно, в силу  теоремы \ref{Tf2} отображение $F$ имеет неподвижную точку. Можно проверить, что неподвижными точками будут все частичные пределы последовательности $(a_i)_\NA$. Такие пределы, конечно, существуют в силу компактности $X$. 

В связи с рассмотренным примером отметим также, что

--- отображение $F$ не является $k$-сжимающим --- достаточно рассмотреть переход аргумента через центр какого ни будь из отрезков $[a_i,b_i]$ --- и поэтому не применима теорема Надлера \cite{Nadler1969};

--- отображение $F$ не является $\alpha$-накрывающим, так как не является сюръективным, и поэтому не применима теорема о совпадающих точках \cite{Arut2007};

--- значения отображения $F$ не выпуклы и могут принимать значение $\varnothing$, поэтому не применима теорема Какутани \cite{Kakutani1941}.

В отношении теоремы Какутани \cite{Kakutani1941} можно заметить, что условия второй части теоремы \ref{Tf2} слабее условий этой теоремы и ее обобщений (см., например, \cite{Glicks1952}). В самом деле, с одной стороны, так как условие второй части является необходимым для существования неподвижной точки отображения $F$, при выполнении условий теоремы Какутани (как и любых других достаточных условий существования неподвижной точки отображения $F$) будут выполняться условия второй части теоремы \ref{Tf2}. Значит, эти условия не сильнее условий теоремы Какутани. С другой стороны, пример \ref{exaTf2} показывает, что условия второй части теоремы \ref{Tf2} не совпадают с условиями теоремы Какутани.

%%%%%%%%%%%%%%%%%%%%%%%%%%%%%%%%%%%%%%%%%%%%%%%%%%%%%%%%%%%%%
\section{Неподвижные точки и итерации функций}

Свойство идемпотентности некоторого оператора (функции) $F$, действующего из $X$ в $X$, 
$$
F(x)=F(F(x))\qquad \forall x\in X
$$
дает множество неподвижных точек этого оператора $\fix{F}=\{F(x)\mid x\in X\}$.
Вообще, исходя из содержательного смысла понятий идемпотентности и неподвижной точки, можно сказать, что неподвижные точки оператора это его область идемпотентности.

Пусть для некоторого подмножества порядковых чисел определены  $\alpha$--итерации $X\ni x\mapsto F^\alpha(x)\in X$ (конечные или бесконечные) оператора $F$. Скажем, что итерации вырождаются при $\alpha$, если 
$$
F^\beta(x)=F^{\alpha}(x)\qquad \forall\beta\ge\alpha, \forall x\in X.
$$
Этот случай также дает описание множества неподвижных точек функции $F$:
\beq\label{fix-F-alfa}
\fix{F}=\{F^\alpha(x)\mid x\in X\}.
\eeq
Часто идемпотентен не сам оператор, а некоторая его итерация:
$$
F^\alpha(x)=F^\alpha(F^\alpha(x))\qquad \forall x\in X.
$$
В этом случае, мы получаем описание лишь множества неподвижных точек оператора $F^\alpha$:
$$
\fix{F^\alpha}=\{F^\alpha(x)\mid x\in X\}.
$$
Однако, в случае бесконечного $\alpha$, мы, неформально говоря, имеем $\alpha+1=\alpha$ и при некоторых дополнительных условиях вновь можем получить представление множества неподвижных точек $F$ в виде \fref{fix-F-alfa}.

В теории обобщенной выпуклости \cite{Soltan1984}, такой результат получается для оператора (выпуклой) предоболочки при наличии в $X$ частичного порядка, а также свойств изотонности и экстенсивности этого оператора:
$$
((x\le y)\myimp(F(x)\le F(y)))\myand(x\le F(x))\qquad\myll x,y\in X.
$$
Заметим, что содержательно этот результат перекликается с певой частью --- существованием --- теоремы Тарского о неподвижной точке \cite{Tarski1955}. Более того, в работе \cite{Eche2005} приводится доказательство этой части теоремы Тарского, использующее технику аналогичную технике доказательства \cite[Теорема 1.4]{Soltan1984}.

Примерами описанных выше ситуаций с идемпотентными операторами могут послужить операторы замыкания \cite[гл.5]{Birkhoff1984} или выпуклой оболочки \cite{Soltan1984}. Случаи вырождения итераций исследованы в различных вариантах метода программных итераций (см., например, \cite{SubChe81, ChenDyat87} и библиографию в этих работах). 

Далее мы приведем результаты о способах <<построения>> неподвижных точек, примыкающие к описанным методам, но не сводящиеся к ним.

Следуя \cite{Engelking1986} обозначим $(X,\lambda)$, где 
$$
L\ni (x_i)_\NA\mapsto\lambda(x_i)_\NA\in X,\qquad L\in 2^{X^\NA}\setminus\varnothing,
$$
пространство со сходимостью $\lambda$ --- оператором, приписывающим отдельным последовательностям $(x_i)_\NA$ из $X$ значения $\lambda(x_i)_\NA$ их пределов.

Пусть сходимость  $\lambda$ удовлетворяет условиям (L1)---(L4): $\myll x\in X$, $ \myll (x_i)_\NA,(y^i_{k})_\NA\in L$
$$
(\myll i\in\NA\ x_i=x)\myimp(\lambda(x_i)_\NA=x)  \eqno {\textup{(L1)}}
$$
$$
(\lambda(x_i)_\NA=x)\myimp((\myll i_1<i_2<\ldots)\lambda(x_{i_k})_\NA=x) \eqno {\textup{(L2)}}
$$
$$
(\lambda(x_i)_\NA\neq x)\myimp((\myxst i_1<i_2<\ldots)(\myll k_1<k_2<\ldots) \lambda(x_{i_{k_n}})_\NA\neq x) \eqno {\textup{(L3)}}
$$
$$
(\lambda(x_i)_\NA= x)\myand(\lambda(y^i_{k})_\NA= x_i)\myimp((\myxst i_1,i_2,\ldots)(\myxst k_1,k_2,\ldots) \lambda(y^{i_n}_{k_n})_\NA= x) \eqno {\textup{(L4)}}
$$

Тогда \cite[п.1.7.18, стр.108]{Engelking1986} существует $T_1$ топология на $X$, обозначим ее $\tau_\lambda$, в которой сходящимися будут только последовательности из $L$ и для всякой такой последовательности $(x_i)_\NA\in L$ будет выполняться равенство
$$
\lim_{i\to \infty}x_i=\lambda(x_i)_\NA,
$$
в котором предел понимается как предел последовательности в топологическом пространстве $(X,\tau_\lambda)$.

Для всякого $n\in\NA$ произвольной функции $X\ni x\mapsto f(x)\in{X}$ определим $n$--ную итерацию $X\ni x\mapsto f^n(x)\in{X}$ функции $f$ по индукции:
$$
f^0(x)\mydef x, \qquad f^{n}(x)\mydef f(f^{n-1})(x)\qquad \forall x\in X.
$$

\begin{teoRU}\label{teo-limit-eq-fix}
Пусть задана функция  $X\ni x\mapsto f(x)\in{X}$ и пространство со сходимостью $(X,\lambda)$, где 
$$
L\ni (x_i)_\NA\mapsto\lambda(x_i)_\NA\in X,\qquad L\in 2^{X^\NA}\setminus\varnothing.
$$
Пусть сходимость $\lambda$ удовлетворяет условиям  (L1)---(L4), функция $f$ непрерывна в топологии $\tau_\lambda$, порожденной этой сходимостью, а последовательности, образованные итерациями функции $f$ содержатся в области определения сходимости $\lambda$:
\beq\label{it-in-lim}
(f^n(x))_\NA\in L\qquad\forall x\in X.
\eeq 
Тогда выполняется равенство
\beq\label{limit-eq-fix}
\fix{f}=\{\lambda(f^n(x))_\NA\mid x\in X\}.
\eeq
В частности, $\fix{f}\neq\varnothing$.
\end{teoRU}

Пусть $(X,\myle)$ --- полное частично упорядоченное множество (п.ч.у.м.) \cite[п.1.2.1]{BAREN1985R}, то есть частично упорядоченное множество, в котором существует  минимальный элемент $\bot\mydef\sup\varnothing$, и в котором каждое направленное подмножество $D\subset X$ имеет точную верхнюю грань $\sup D$.

Напомним, что топология Скотта  \cite[п.1.2.3]{BAREN1985R} (обозначим ее $\tau_S$) определяется на п.ч.у.м.  $(X,\myle)$ следующим образом:   множество $O\subseteq X$ является открытым, если и только если для любых $x,y\in X$ и направленного множества $D\in2^X$
$$
(x\in O)\myand(x\myle y)\myimp(y\in O),
$$
$$
(\sup D\in O)\myimp(D\cap O\neq\varnothing).
$$

\begin{teoRU}\label{teo-limit-eq-fix-2}
Пусть задана функция  $X\ni x\mapsto f(x)\in{X}$ и функция $f$ непрерывна в топологии $\tau_S$, определенной на п.ч.у.м.  $(X,\myle)$. 
Тогда выполняется равенство
\beq\label{limit-eq-fix-2}
\fix{f}=\{\sup \{f^n(x)\mid n\in\NA\}\mid x\myle f(x)\}.
\eeq
В частности, $f$ имеет наименьшую неподвижную точку $\sup \{f^n(\bot)\mid n\in\NA\}$.
\end{teoRU}

\begin{corRU}\label{cor-limit-eq-fix-2}
Пусть функция $X\ni x\mapsto f(x)\in{X}$ непрерывна в топологии $\tau_S$ и экстенсивна:
$$
x\myle f(x)\qquad x\in X.
$$ 
Тогда выполняется равенство
\beq\label{limit-eq-fix-3}
\fix{f}=\{\sup \{f^n(x)\mid n\in\NA\}\mid x\in X\}.
\eeq
В частности, $f$ имеет наименьшую неподвижную точку $\sup\{f^n(\bot)\mid n\in\NA\}$.
\end{corRU}

\begin{remRU}\label{rem-limit-eq-fix-2}
Следствие \ref{cor-limit-eq-fix-2} перекликается с результатами \cite[Теоремы 2.1, 2.2]{ChenDyat87}, сформулированными для частично упорядоченного множества $(\mathcal X,\subseteq)$, где $\mathcal X\subseteq 2^Y$ --- некоторое семейство подмножеств $Y$, замкнутое относительно монотонной сходимости множеств. 

Заметим также, что монотонная сходимость множеств, используемая в результатах \cite[Теоремы 2.1, 2.2]{ChenDyat87}, не удовлетворяет условиям теорем \ref{teo-limit-eq-fix}, \ref{teo-limit-eq-fix-2} (см. пример \ref{ex-seq-not-L4}).
\end{remRU}

\begin{exaRU}\label{ex-seq-not-L4}
Пусть $X\mydef 2^{\RA^2}$,  частичный порядок задан отношением включения:
$$
(x\myle y)\myeqv (x\supseteq y),\qquad \myll x,y\in X,
$$
и сходящимися являются только цепи:
$$
((x_i)_\NA\in L)\myeqv (x_i\myle x_{i+1}, \myll i\in\NA), \qquad \lambda(x_i)_\NA\mydef\cap_{i\in\NA}x_i.
$$
В качестве $\mathcal X$ рассмотрим множества замкнутые в топологии, порожденной евклидовой нормой. Понятно, что это семейство инвариантно относительно указанной сходимости. 

Пусть множества $x_i, y_k^i\in \mathcal X$, $i,k\in\NA$ заданы следующим образом
$$
x_i\mydef\{(\alpha,0)\mid \alpha\in[0,1/i]\}, \quad   y_k^i\mydef x_i\cap\{(1/i,\beta)\mid \beta\in[0,1/k]\}, \qquad i,k\in\NA.
$$
Тогда выполняется посылка условия (L4). Однако, при выборе любой цепи $(z_i)_\NA$ из множества $\{y_k^i\mid i,k\in\NA\}$, для некоторого $n\in\NA$ верны соотношения
$$
\lambda(z_i)_\NA\supseteq\{(\alpha,0)\mid \alpha\in[0,1/n]\}\neq\{(0,0)\}=\lambda(x_i)_\NA.
$$

\end{exaRU}

\section{Приложение к теоремам о равновесиях}

Используя известный прием, связывающий положения равновесия с неподвижными точками многозначных отображений (см. утверждения \ref{lem_fix_eq_saddle}, \ref{lem_fix_eq_NE}), из приведенных теорем о неподвижных точках можно получить новые варианты теорем о равновесиях.

1. Так из теоермы \ref{TfO} можно получить еще одну теорему о критерии равенства минимакса и максимина в условиях сходных с условиями теоремы Фана  \cite{FanKy1953}.

Пусть  $\UB$, $\VB$ --- (топологические) пространства чистых стратегий двух игроков, на произведении  
$$
X\mydef\UB\times\VB
$$ 
заданы стандартная топология произведения и функция исходов $\varphi: \UB\times\VB\mapsto\RA$ со скалярными значениями. Игрок, выбирающий стратегии $u\in\UB$, стремится минимизировать значение исхода игры, игрок, выбирающий $v\in\VB$, стремится к увеличению исхода. Обозначим $S(\varphi)\in2^\UB\times2^\VB$ множество седловых точек функции $\varphi$, то есть множество пар $(u_*,v_*)\in\UB\times\VB$, удовлетворяющих условиям
\beq\label{def-sed}
\varphi(u_*,v)\le\varphi(u,v_*)\qquad\myll (u,v)\in\UB\times\VB.
\eeq 
Введем в рассмотрение многозначное отображение $X\ni(u,v)\mapsto \Fg{F}{\varphi}(u,v)\in2^{X}$ вида
\beq\label{F-for-fi}
\Fg{F}{\varphi}(u,v)\mydef\argmin_{u'\in\UB}\varphi(u',v)\times\argmax_{v'\in\VB}\varphi(u,v'),\qquad(u,v)\in X.
\eeq
\begin{teoRU}\label{teoFC-TfKK}
Пусть $\UB$, $\VB$ --- компактные Хаусдорфовы пространства; для всякого $v\in\VB$ функция  $\varphi(\cdot,v):\UB\mapsto\RA$ полунепрерывна снизу на $\UB$; для всякого $u\in\UB$ функция  $\varphi(u,\cdot):\VB\mapsto\RA$ полунепрерывна сверху на $\VB$. При этих условиях для множество $S(\varphi)$ седловых точек функции $\varphi$ выполняются равенства 
\beq\label{sed-eq-fix-eq-nocycle-o}
S(\varphi)=\fix{\Fg{F}{\varphi}}=\bigcap_{(O_\iota)_{ I}\in \Ofo} \bigcup_{\iota\in I}\Fg{\barf{F}}{\varphi}(O_\iota)=\bigcap_{(O_\iota)_{ I}\in \Ofc} \bigcup_{\iota\in I}\Fg{\barf{F}}{\varphi}(O_\iota).
\eeq
В частности, равенство 
\beq\label{s1}
\max_{v\in\VB}\min_{u\in\UB}\varphi(u,v)=\min_{u\in\UB}\max_{v\in\VB}\varphi(u,v)
\eeq
выполнено тогда и только тогда, когда в произвольном покрытии $(O_\iota)_I\in\Ofc$ найдется множество $O_{\bar\iota}$, $\bar\iota\in I$, содержащее два последовательных приближения Курно
$$
x,x'\in O_{\bar\iota},\qquad x'\in\Fg{F}{\varphi}(x).
$$
\end{teoRU}

%%
%%\begin{remRU}
%%Для доказательства теоремы \ref{teoFC-TfKK} достаточно установить, что условие \fref{m1} эквивалентно условию \fref{con-cov-not-empt} для отображения \fref{F-for-fi}. Будет приведено формальное доказательство этого факта, использующее теоремы \ref{teoFC} и \ref{TfKK}. Если будет найдено непосредственное доказательство эквивалентности этих условий, то тем самым будет найдено еще одно доказательство теоремы Фан Цзы.
%%\end{remRU}

2. Переходя к более общему случаю игры получим еще один вариант теоремы о существовании равновесия по Нэшу. Введем необходимые обозначения. 

Пусть $(X,J)$ игра с $n$ игроками в нормальной форме:
$$
X\mydef X_1\times\ldots\times X_n,\quad J\mydef (J_1,\ldots J_n),
$$
$$
(y,x_{-i})\mydef(x_1,\ldots,x_{i-1},y,x_{i+1},\ldots,x_n),\quad y\in X_i,\ x=(x_1,\ldots,x_n)\in X.
$$
Здесь $J_i$ --- функция выигрыша $i$-го игрока: $J_i:X\mapsto\RA$, $i\in\nint{1}{n}$ и в случае, когда $X_i$ --- топологические пространства, полагаем, что на произведении $X$ действует стандартная топология произведения.
Будем обозначать $N(J)$ множество элементов $x^*\in X$, в которых достигается равновесие по Нэшу, то есть выполнены условия
\beq\label{NE}
J_i(y_i,x^*_{-i})\le J_i(x^*),\qquad\myll y_i\in X_i, \myll i\in\nint{1}{n}.
\eeq
По аналогии с \fref{F-for-fi} введем многозначное отображение $X\ni x\mapsto\Fg{F}{(X,J)}(x)\in 2^X$:
\beq\label{F-for-J}
\Fg{F}{(X,J)}(x)\mydef\argmax_{y_1\in X_1}J_1(y_1,x_{-1})\times\ldots\times\argmax_{y_n\in X_n}J_n(y_n,x_{-n}),\qquad x\in X.
\eeq

\begin{teoRU}\label{NEQ-TfKK}
Пусть $X_i$, $i\in\nint{1}{n}$ --- компактные Хаусдорфовы пространства; для всех $x\in X$, $i\in\nint{1}{n}$ функция  $J_i(\cdot,x_{-i}):X_i\mapsto\RA$ полунепрерывна сверху на $X_i$. При этих условиях для множества $N(J)$ Нэшевских рановесий выполняются соотношения
\beq\label{nesh-eq-fix-eq-nocycle-o}
N(J)=\fix{\Fg{F}{(X,J)}}=\bigcap_{(O_\iota)_{ I}\in \Ofo} \bigcup_{\iota\in I}\Fg{\barf{F}}{(X,J)}(O_\iota)=\bigcap_{(O_\iota)_{ I}\in \Ofc} \bigcup_{\iota\in I}\Fg{\barf{F}}{(X,J)}(O_\iota).
\eeq
В частности, равновесие по Нэшу \fref{NE} достигается тогда и только тогда, когда в произвольном покрытии $(O_\iota)_I\in\Ofc$ найдется множество $O_{\bar\iota}\in (O_\iota)_I$, содержащее два последовательных приближения Курно
\beq\label{kur-in-O}
(\myll (O_\iota)_I\in\Ofc)(\myxst \bar\iota\in I)(\myxst x,x'\in O_{\bar\iota})\qquad x'\in\Fg{F}{(X,J)}(x).
\eeq
\end{teoRU}

\begin{remRU}
В случае, когда топология пространства $X$ метризуема, условие \fref{kur-in-O} естественным образом трансформируется в условие \fref{x-close-to-Fx}, принимая вид
$$
(\myll \delta>0)(\myxst x_\delta\in X)\qquad  \dis(x_\delta,\Fg{F}{(X,J)}(x_\delta))\le\delta.
$$
\end{remRU}

\begin{remRU}
Отметим, что теорема \ref{NEQ-TfKK} также, как теорема \ref{teoFC-TfKK} и теорема Фана \cite{FanKy1953} приобрела характер критерия.
\end{remRU}

\section{Доказательства}

\subsection{Доказательство теоремы \ref{LOx}}
Включение 
\beq\label{LOx-inc-0}
\fix{F}\subset\bigcap_{(O_\iota)_{ I}\in \Ox} \bigcup_{\iota\in I}\barf{F}(O_\iota).
\eeq
следует из \fref{Lf2-inc2}. В самом деле, для любого $(O_\iota)_{ I}\in \Ox$ в силу \fref{Lf2-inc2} имеем равенства
$$
\fix{F}=\fix{F}\cap \bigcup_{\iota\in I}O_\iota=\fix{F}\cap \bigcup_{\iota\in I}\barf{F}(O_\iota),
$$
и, следовательно, 
$
\fix{F}\subset \bigcup_{\iota\in I}\barf{F}(O_\iota).
$

Покажем обратное включение. Так как 
$$
\{\fix{F},\{x\}\mid x\in X, x\not\in F(x)\}\in \Ox,
$$
из принадлежности $\bar x$ правой части \fref{fix-eq-nocycle-x} получим, в частности,  включение
$$
\bar x\in \barf{F}(\fix{F})\cup\bigcup_{x\in X\setminus\fix{F}}\barf{F}(\{x\}),
$$
из которого с учетом \fref{Lf1-inc1} и равенств
$
\barf{F}(\{x\})=\varnothing$ при $x\in X\setminus\fix{F}
$
следуют включения
$
\bar x\in \barf{F}(\fix{F})\subset\fix{F}.
$

\subsection{Доказательство теоремы \ref{TfO} и следствия \ref{TfTS}}

%tag proof 2.2

В силу утверждения \ref{LOx} и включений $\Ofo,\Ofc\subset\Ox$ для доказательства утверждения достаточно установить, что 
\beq\label{incl1}
\fix{F}\supset\bigcap_{(O_\iota)_{ I}\in \Ofo} \bigcup_{\iota\in I}\barf{F}(O_\iota),
\eeq
\beq\label{incl2}
\fix{F}\supset\bigcap_{(O_\iota)_{ I}\in \Ofc} \bigcup_{\iota\in I}\barf{F}(O_\iota).
\eeq

1. Обратимся к обоснованию \fref{incl1}. Пусть $\bar x$ принадлежит правой части \fref{incl1}.
Для произвольной локальной базы топологии $(O_s(\bar x))_{S}$ компакта $X$ в точке $\bar x$ построим семейство  покрытий $(O_{s\iota})_{ I_s}$, $s\in S$ из $\Ofo$ такое, что для всякого $s\in S$ выполнены соотношения
\beq\label{TfH-eq0}
O_s(\bar x)=O_{s\iota_s},\quad \iota_s\in I_s,\qquad \bar x\not\in O_{s\iota}, \quad \iota\in I_s\setminus\{\iota_s\}.
\eeq
Выберем произвольное $s\in S$ и, пользуясь свойством Хаусдорфа, для каждого элемента $y\in X\setminus O_s(\bar x)$ выберем открытые непересекающиеся окрестности $O_{sy}$ и $O_{s\bar x}$ элементов $y$ и $\bar x$. Семейство 
\beq\label{TfH-eq2}
\{O_s(\bar x), O_{sy}\mid y\in X\setminus O_s(\bar x)\}
\eeq
является открытым покрытием $X$. В качестве покрытия $(O_{s\iota})_{ I_s}$ при выбранном $s$ возьмем конечное подпокрытие покрытия \fref{TfH-eq2}. По построению оно обладает всеми заявленными свойствами.

Из соотношений \fref{Lf2-inc1}, \fref{TfH-eq0} и предположения о принадлежности  $\bar x$ правой части \fref{fix-eq-nocycle-o} для всех $s\in S$ следуют включения
\beq\label{TfH-eq1}
\bar x
\in \bigcup_{\iota\in I_s}\barf{F}(O_{s\iota})\setminus\bigcup_{\iota\in I_s\atop\iota\neq\iota_s}O_{s\iota}
\subset  \bigcup_{\iota\in I_s}\barf{F}(O_{s\iota})\setminus\bigcup_{\iota\in I_s\atop\iota\neq\iota_s}\barf{F}(O_{s\iota})
\subset \barf{F}(O_s(\bar x))=\Biggl(\bigcup_{y\in O_s(\bar x)}F(y)\Biggl) \bigcap O_s(\bar x).
\eeq

Покажем, что $\bar x \in\fix{F}$. 
Предположим противное: $\bar x\notin F(\bar x)$. Из замкнутости графика $F$ следует зам\-к\-нутость множества $F(\bar x)\subset X$. Отсюда в силу компактности и хаусдорфовости $X$ следует существование открытых окрестностей $O'(\bar x)$,  $O'(F(\bar x))$ точки $\bar x$ и замкнутого множества $F(\bar x)$ таких, что
\beq\label{TfO-eq1}
O'(\bar x)\bigcap O'(F(\bar x))=\varnothing.
\eeq
В силу замкнутости графика отображения $F$, компактности и хаусдорфовости пространства  $X$ отображение $F$ является полунепрерывным сверху (см.  \cite[\S43,I,Теорема 4]{Kura1969}, \cite[Теорема 1.2.32]{BorGelMysObu2011}), то есть для окрестности   $O'(F(\bar x))$ найдется окрестность  $O''(\bar x)$ такая, что при всех $y\in O''(\bar x)$ будут выполняться включения
\beq\label{TfO-eq2}
F(y)\subset O'(F(\bar x)).
\eeq
Из определения семейства $(O_s(\bar x))_{s\in S}$ следует, что для некоторого $\bar s\in S$ будет выполнено
\beq\label{TfO-eq3}
O_{\bar s}(\bar x)\subset O'(\bar x)\cap O''(\bar x).
\eeq
Cоотношения \fref{TfO-eq1}, \fref{TfO-eq2}, \fref{TfO-eq3} дают включения
$$
\Biggl(\bigcup_{y\in O_{\bar s}(\bar x)}F(y)\Biggl) \bigcap O_{\bar s}(\bar x)
\subset\Biggl(\bigcup_{y\in O''(\bar x)}F(y)\Biggl) \bigcap O'(\bar x)
\subset O'(F(\bar x)) \bigcap O'(\bar x)  =\varnothing,
$$
которые  ведут к противоречию с \fref{TfH-eq1} при $s=\bar s$.

2. Для обоснования \fref{incl2} докажем два вспомогательных утверждения.

\begin{lemRU}\label{Lf3}
Пусть $X$ --- компактное хаусдорфово пространство. Для любого элемента $x\in X$ и любой окрестности $O\in\tau(X)$, $x\in O$ найдется замкнутое множество $B$ и окрестность $O'\in\tau(X)$ такие, что $x\in  O'\subset B\subset O$.
\end{lemRU}

\bpr
Пусть $x\in O\in\tau(X)$. Так как $X\setminus O$ замкнуто и не содержит $x$ в силу компактности и хаусдорфовости $X$ найдутся $O',O''\in\tau(X)$ удовлетворяющие
$
x\in O'$, $X\setminus O\subset O''$, $O'\cap O''=\varnothing.
$

Тогда для замкнутого множества $B\mydef X\setminus O''$ и окрестности $O'$ выполняются искомые соотношения: $x\in O'\subset B\subset O$.
\epr

\begin{lemRU}\label{Lf4}
Пусть $X$ --- компактное хаусдорфово пространство. Для любых $x\in X$, $O\in\tau(X)$ и замкнутого множества $B\in X$, удовлетворяющих включениям $x\in O\subset B$. Тогда найдется покрытие $(O_\iota)_I\in \Ofc$ такое, что для некоторого $\bar\iota\in I$ выполнены соотношения
\beq\label{Lf4-eq1}
B=O_{\bar\iota},\qquad x\in X\setminus\bigcup_{\iota\in I\setminus\{\bar\iota\}}O_\iota.
\eeq
\end{lemRU}

\bpr
Пусть $x\in X$, $O, B\in2^X$ удовлетворяют условиям этой леммы.  Для любой точки $y\in X\setminus\{x\}$, пользуясь свойством хаусдорфа пространства $X$, построим окрестности $O_{xy},O_{yx}\in\tau(X)$ такие, что 
$
x\in O_{xy}\subset O$, $y\in O_{yx}$, $O_{xy}\cap O_{yx}=\varnothing.
$

Пользуясь компактностью $X$ из открытого покрытия $\{O_{xy},O_{yx}\mid y\in X\setminus\{x\}\}$ выделим конечное подпокрытие
$
\{O_{xy_i}, O_{y_jx}\mid i\in I, j\in J\}
$
и рассмотрим конечное семейство замкнутых множеств
\beq\label{Lf4-eq2}
\{B,B_j\mid j\in J\},
\eeq
здесь $B_j$ обозначает замыкание множества $O_{y_jx}$. 
Так как $O_{xy_i}\subset O\subset B$, $i\in I$ это семейство есть покрытие $X$. По построению 
$$
\left(\bigcup_{j\in J}O_{y_jx}\right)\bigcap\left(\bigcap_{i\in I}O_{xy_i}\right)=\varnothing.
$$
Отсюда, учитывая, что $x\in \bigcap_{i\in I}O_{xy_i}$, получим
$
x\in X\setminus\bigcup_{j\in J}B_j.
$
Мы проверили последнее требуемое свойство покрытия \fref{Lf4-eq2}.
\epr

3. Докажем включение \fref{incl2} <<от противного>>: пусть $x$ принадлежит правой части \fref{incl2} и $x\not\in F(x)$. В силу замкнутости графика $F$ множество $F(x)$ замкнуто. Значит в силу компактности и хаусдорфовости имеются непересекающиеся окрестности $O'(x),O'(F(x))\in \tau(X)$ точки $x$ и множества $F(x)$:
$
O'(x)\cap O'(F(x))=\varnothing.
$
В силу замкнутости графика $F$ и условий на $X$ отображение $F$ является полунепрерывным сверху по включению (см.  \cite[\S43,I,Теорема 4]{Kura1969}, \cite[Теорема 1.2.32]{BorGelMysObu2011}), то есть существует окрестность $O''(x)\in\tau(X)$ такая, что
$
\bigcup_{y\in O''(x)}F(y)\subset O'(F(x)).
$

Пусть в силу леммы \ref{Lf3} найдены открытая окрестность $O$ и замкнутое множество $B$, удовлетворяющие 
$
x\in O\subset B\subset O'(x)\cap O''(x).
$
Пусть в силу леммы \ref{Lf4} покрытие  $(O_\iota)_I\in \Ofc$ удовлетворяет условиям \fref{Lf4-eq1}. Тогда из этих условий, предположения 
$
x\in\bigcap_{(O_\iota)_{ I}\in \Ofc} \bigcup_{\iota\in I}\barf{F}(O_\iota)
$
и свойства \fref{Lf2-inc1} следует
$$
x\in\left(\bigcup_{\iota\in I}\barf{F}(O_\iota)\right)\setminus\left(\bigcup_{\iota\in I\setminus\{\bar\iota\}}O_\iota\right)
\subset\left(\bigcup_{\iota\in I}\barf{F}(O_\iota)\right)\setminus\left(\bigcup_{\iota\in I\setminus\{\bar\iota\}}\barf{F}(O_\iota)\right)
\subset\barf{F}(O_{\bar\iota})=\barf{F}(B).
$$
С другой стороны, из построения $B$ имеем
$$
\barf{F}(B)\mydef B\bigcap \left(\bigcup_{y\in B} F(y)\right)
\subset O'(x)\bigcap \left(\bigcup_{y\in O''(x)} F(y)\right)
\subset O'(x)\bigcap O'(F(x))=\varnothing.
$$
Получили противоречие: $x\in\varnothing$. Первая часть теоремы доказана.

4. Для проверки второй части теоремы сначала установим, что при выполнении условия \fref{con-cov-not-empt} множества 
\beq\label{TfKK-eq1}
\bigcup_{\iota\in I}\barf{F}(O_\iota),\qquad (O_\iota)_I\in \Ofc
\eeq
замкнуты и центрированы. Из этого свойства в совокупности с компактностью множества $X$ будет следовать непустота правой части \fref{incl2}.

Замкнутость этих множеств следует из определений семейства $\Ofc$, отбражения $\barf{F}$ и замкнутости отображения $F$ (см. \cite[Теорема 1.2.33]{BorGelMysObu2011}): для любого замкнутого множества $Y\subset X$ множество $\cup_{y\in Y}F(y)$ замкнуто.

Проверим центрированность множеств \fref{TfKK-eq1}: пусть $(O_{k\iota})_{I_k}\in \Ofc$, $k\in\nint{1}{n}$, $n\in\NA$. Положим 
$$
(\bar O_{\bar\iota})_{\bar I}\mydef\{O_{1\iota_1}\cap\ldots\cap O_{n\iota_n}\mid \iota_k\in I_k, k\in\nint{1}{n}\}.
$$
Иначе говоря $(\bar O_{\bar\iota})_{\bar I}$ есть точная нижняя граница конечного множества покрытий $(O_{k\iota})_{I_k}$ по отношению вписанности.

По построению $(\bar O_{\bar\iota})_{\bar I}\in \Ofc$ и вписано во все разбиения $(O_{k\iota})_{I_k}$:
$
(\bar O_{\bar\iota})_{\bar I}\poc(O_{k\iota})_{I_k}$, $k\in\nint{1}{n}. 
$
Cледовательно (см. \fref{Lf2-inc4}), при любом $k\in\nint{1}{n}$ выполнено включение
$
\bigcup_{\bar\iota\in\bar I}\barf{F}(\bar O_{\bar\iota})\subset \bigcup_{\iota\in I_k}\barf{F}(O_{k\iota}).
$
Отсюда, с учетом условия $\bigcup_{\bar\iota\in\bar I}\barf{F}(\bar O_\iota)\neq\varnothing$, получим
$
\varnothing\neq\bigcup_{\bar\iota\in\bar I}\barf{F}(\bar O_{\bar\iota})\subset \bigcap_{k\in\nint{1}{n}}\bigcup_{\iota\in I_k}\barf{F}(O_{k\iota}).
$

Таким образом, подмножества \fref{TfKK-eq1} компактного пространства $X$ замкнуты и центрированы. Следовательно, 
$
\fix{F}=\bigcap_{(O_\iota)_{ I}\in \Ofc} \bigcup_{\iota\in I}\barf{F}(O_\iota)\neq\varnothing.
$

Обратно, пусть $\bar x\in\fix{F}$, тогда в любом покрытии $(O_\iota)_I\in\Ox$ найдется множество $O_{\bar\iota}\in(O_\iota)_I$ такое, что $\bar x\in O_{\bar\iota}$. Следовательно, по определению, $\bar x\in\barf{F}(O_{\bar\iota})$ и, значит, $\barf{F}(O_{\bar\iota})\neq\varnothing$. 

Для обоснования следствия \ref{TfTS} заметим, что условие замкнутости графика отображения эквивалентно условию непрерывности функции.

\subsection{Доказательство теоремы \ref{Tf2} и следствия \ref{Tf2S}}

В этом пункте для произвольных $\delta>0$, $x\in X$ и $Y\in2^X$ будем обозначать
$$
O_\delta(x)\mydef\{y\in X\mid \di(x,y)\le\delta\}, \qquad O_\delta(Y)\mydef\{y\in X\mid \dis(y,Y)\le\delta\}.
$$

1. В силу теоремы \ref{LOx} и включений $(O_{k\iota})_{I_k}\subset\Ox$,  $k\in\NA$ для доказательства первой части достаточно установить, что 
$
\fix{F}\supset\bigcap_{k\in\NA} \bigcup_{\iota\in I_k}\barf{F}(O_{k\iota}).
$

Предположим, что $\bar x\in \bigcap_{k\in\NA} \bigcup_{\iota\in I_k}\barf{F}(O_{k\iota})$. Тогда имеется последовательность 
\beq\label{Tf2-eq2}
\bar O_k\mydef O_{k\iota_k}, \quad \bar x\in \barf{F}(\bar O_k)\subset\bar O_k,\quad \diametr(\bar O_{k})\le \delta_k,\qquad k\in\NA.
\eeq
Покажем, что $\bar x \in\fix{F}$. 
Предположим противное: $\bar x\notin F(\bar x)$. Положим 
$
\bar\varepsilon\mydef\dis(\bar x, F(\bar x))/3.
$
В силу замкнутости графика множество $F(\bar x)$ замкнуто и поэтому $\bar\varepsilon>0$. Выберем $\bar\delta\in(0,\bar\varepsilon]$ так, чтобы для всех $x\in O_{\bar\delta}(\bar x)$ выполнялись соотношения
$
F(x)\subset O_{\bar\varepsilon}(F(\bar x)).
$
Это можно сделать, так как из замкнутости графика следует полунепрерывность сверху по включению отображения $F$  (см.  \cite[\S43,I,Теорема 4]{Kura1969}, \cite[Теорема 1.2.32]{BorGelMysObu2011}). При этом из выбора величин $\bar\delta$, $\bar\varepsilon$ следует соотношение
$
O_{\bar\delta}(\bar x)\bigcap O_{\bar\varepsilon}(F(\bar x))=\varnothing.
$
Из этих построений получим 
\beq\label{Tf2-eq1}
\barf{F}(O_{\bar\delta}(\bar x))\mydef O_{\bar\delta}(\bar x)\bigcap\left(\bigcup_{x\in O_{\bar\delta}(\bar x)}F(x)\right)\subset O_{\bar\delta}(\bar x)\bigcap O_{\bar\varepsilon}(F(\bar x))=\varnothing.
\eeq
Выберем $\bar k\in\NA$ таким, чтобы выполнялось включение $\bar O_{\bar k}\subset O_{\bar\delta}(\bar x)$. Это можно сделать в силу \fref{dim-diam-cov}. Из этих включений и соотношений \fref{Tf2-eq1}, \fref{Lf1-inc2} получим
$
\barf{F}(\bar O_{\bar k})\subset \barf{F}(O_{\bar\delta}(\bar x))=\varnothing,
$
что противоречит \fref{Tf2-eq2} при $k=\bar k$. Первая часть теоремы доказана.

2.  Обратимся ко второй части утверждения. Необходимость условия \fref{x-close-to-fx} следует сразу. Проверим достаточность этого условия. Пользуясь компактностью $X$, построим последовательность $(O_{k\iota})_{I_k}$, $k\in\NA$ открытых вписанных покрытий вида:
\beq\label{mon-cov}
(O_{k+1\iota})_{I_{k+1}}\poc(O_{k\iota})_{I_k},\quad(O_{k\iota})_{I_k}\in \Ofod{1/k} \qquad k\in\NA.
\eeq
Тогда последовательность $(\bar O_{k\iota})_{I_k}$, $k\in\NA$, где $\bar O_{k\iota}$--замыкание в $X$ множества $O_{k\iota}$, замкнутых конечных покрытий будет удовлетворять условиям \fref{dim-diam-cov}:
\beq\label{mon-cov-cl}
(\bar O_{k+1\iota})_{I_{k+1}}\poc(\bar O_{k\iota})_{I_k},\quad(\bar O_{k\iota})_{I_k}\in \Ofcd{1/k} \qquad k\in\NA.
\eeq

В этом случае последовательность замкнутых множеств $\bigcup_{\iota\in I_k}\barf{F}(\bar O_{k\iota})$,  в силу \fref{Lf2-inc4}, будет монотонно убывающей по включению. Пользуясь условием \fref{x-close-to-fx}, проверим непустоту членов этой последовательности и, значит, непустоту их пересечения (по части первой теоремы это есть множество неподвижных точек $F$). 

Для всякого $k\in\NA$ выберем $\delta_k$ равным половине числа Лебега покрытия $(O_{k\iota})_{I_k}$. Тогда из условия  \fref{x-close-to-fx} следует существование пары $x_k,y_k\in X$, такой что
$$
\di(x_k,y_k)\le\delta_k, \quad y_k\in F(x_k),
$$
а из выбора $\delta_k$ --- элемента покрытия $O_{k\bar\iota}\in(O_{k\iota})_{I_k}$ такого, что 
$$
x_k,y_k\in O_{k\bar\iota}\subset \bar O_{k\bar\iota}.
$$
Из последних двух соотношений получим 
$$
\barf{F}(\bar O_{k\bar\iota})\neq\varnothing
$$
и, следовательно,
$$
\bigcup_{\iota\in I_k}\barf{F}(\bar O_{k\iota})\neq\varnothing.
$$

Для обоснования следствия \ref{Tf2S} заметим, что условие замкнутости графика отображения эквивалентно условию непрерывности функции.

\subsection{Доказательство теоремы \ref{teo-limit-eq-fix}}

Пусть $x\in\fix{f}$. Тогда $f^n(x)=x$ при всех $n\in\NA$. Отсюда, в силу (L1) имеем $\lambda(f^n(x))_\NA=x$, то есть $x$ лежит в правой части \fref{limit-eq-fix}.

Обозначим $f^\infty(x)\mydef\lambda(f^n(x))_\NA$ и пусть теперь $x\in\{f^\infty(x)\mid x\in X\}$. Значит имеется $ x'\in X$ такой, что для последовательности $(x_i)_\NA$:
$$
x_i\mydef f^i(x')\qquad\forall i\in\NA,
$$
выполнено 
\beq\label{llimxi}
\lambda(x_i)_\NA=x
\eeq 
и, следовательно, в силу определения топологии $\tau_\lambda$, выполняется равенство
\beq\label{limxi}
\lim_{i\to\infty}x_i=x.
\eeq
Обозначим $y\mydef f(x)$ и выберем произвольную открытую окрестность $O_y\in\tau_\lambda$ точки $y$. В силу непрерывности функции $f$ в топологии $\tau_\lambda$ прообраз $O_x\mydef f^{-1}(O_y)$ этой окрестности при отображении $f$ есть открытая окрестность точки $x$: $x\in O_x\in \tau_\lambda$. Следователно, в силу \fref{limxi}, найдется $N(O_x)\in\NA$ такое, что 
$$
x_n\in O_x\qquad \forall n\ge N(O_x).
$$
Тогда для $N(O_y)\in\NA$, равного $N(O_y)\mydef N(O_x)+1$, имеем включения 
$$
x_n\in O_y\qquad \forall n\ge N(O_y).
$$
Таким образом, $\lim_{i\to\infty}x_i=y$. В силу определения топологии $\tau_\lambda$, это равенство выполняется тогда и только тогда, когда
\beq\label{llimxiy}
\lambda(x_i)_\NA=y,
\eeq
Из \fref{llimxiy}, \fref{llimxi} и однозначности определения $\lambda(x_i)_\NA$ получаем равенство $f(x)\mydef y=x$. 
Следовательно, $x\in\fix{f}$.

Так как множество $X$ предполагается не пустым, правая часть \fref{limit-eq-fix} в силу \fref{it-in-lim} также является не пустым множеством. Следовательно, $\fix{f}\neq\varnothing$.

\subsection{Доказательства теоремы \ref{teo-limit-eq-fix-2} и следствия \ref{cor-limit-eq-fix-2}}

Пусть $x\in\fix{f}$. Тогда $(f^n(x))_\NA=\{x\}$. Синглетон $\{x\}$, очевидно, является направленным множеством, $\sup\{x\}=x$ и $x\myle f(x)$. Значит, $x$ лежит в правой части \fref{limit-eq-fix-2}.

Проверим обратное включение. В силу непрерывности функции $f$ в топологии $\tau_S$ (см. \cite[п. 1.2.6, 1.2.7]{BAREN1985R}) эта функция изотонна:
\beq\label{isot-f}
(x\myle y)\myimp(f(x)\myle f(y))\qquad\myll x,y\in X,
\eeq
и для любого направленного множества $D\subset X$ удовлетворяет равенству
\beq\label{sup-cont-f}
\sup\{f(x)\mid x\in D\}=f(\sup D).
\eeq

Зафиксируем произвольный $x\in X$, удовлетворяющий условию $x\myle f(x)$. Заметим, что множество таких элементов не пусто, так как содержит $\bot$. Рассмотрим последовательность $(x_n)_\NA\mydef(f^n(x))_\NA$.  Из определения элемента $x$ и свойства \fref{isot-f} изотонности функции $f$ следует, что последовательность $(x_n)_\NA$ образует цепь:
\beq\label{xn-mono}
x_n\myle x_{n+1}\qquad\myll n\in\NA
\eeq
и уж тем более направленное множество. Значит существует $\bar x\mydef\sup\{x_n\mid n\in\NA\}$. Отметим, что из монотонности \fref{xn-mono} последовательности $(x_n)_\NA$ также следуют равенства 
$$
\bar x=\sup\{x_n\mid n\in\NA,n\ge2\}=\sup\{f(x_n)\mid n\in\NA\}.
$$
Из этих соотношений и \fref{sup-cont-f} получим
$$
\bar x=\sup\{f(x_n)\mid n\in\NA\}=f(\sup\{x_n\mid n\in\NA\})=f(\bar x).
$$
Следовательно, 
$$
\bar x=\sup\{f^n(x)\mid n\in\NA\}\in\fix{f}.
$$

Для обоснования второй части теоремы заметим, что для тех $x\in X$, при которых существует $\sup\{f^n(x)\mid n\in\NA\}$, выполнены неравенства
$$
f^n(\bot) \myle f^n(x)\myle\sup\{f^n(x)\mid n\in\NA\} \qquad\myll n\in\NA.
$$
Тогда по опредлению точной верхней грани $\sup\{f^n(\bot)\mid n\in\NA\}$ имеем 
$$
\sup\{f^n(\bot)\mid n\in\NA\}\myle\sup\{f^n(x)\mid n\in\NA\}.
$$

Доказательство следствия \ref{cor-limit-eq-fix-2} следует доказательству теоремы \ref{teo-limit-eq-fix-2} за исключением понятных изменений в обосновании \fref{xn-mono}.

\subsection{Доказательство теоремы \ref{teoFC-TfKK}}

Докажем вспомогательное утверждение.

\begin{lemRU}\label{lem_fix_eq_saddle} Множество неподвижных точек отображения $\Fg{F}{\varphi}$, то есть множество пар  $(u_*,v_*)\in\UB\times\VB$, удовлетворяющих условию 
\beq\label{FP}
(u_*,v_*)\in \Fg{F}{\varphi}(u_*,v_*),
\eeq
совпадает со множеством седловых точек функции $\varphi(\cdot)$, то есть пар $(u_*,v_*)\in\UB\times\VB$ удовлетворяющих неравенствам 
\beq\label{sedlo}
\varphi(u_*,v)\le\varphi(u,v_*)\qquad\myll(u,v)\in\UB\times\VB.
\eeq
\end{lemRU}

\bpr
Предположим выполнено включение \fref{FP}. Это означает, что справедливы неравенства
\begin{gather}
\varphi(u_*,v_*)\le\myinf_{u'\in\UB}\varphi(u',v_*)\le\varphi(u,v_*),\qquad u\in\UB,\label{FP1}\\
\varphi(u_*,v)\le\sup_{v'\in\VB}\varphi(u_*,v')\le\varphi(u_*,v_*),\qquad v\in\VB.\label{FP2}
\end{gather}
Так как левая часть неравенства \fref{FP1} равняется правой в  \fref{FP2} справедливы соотношения \fref{sedlo}, означающие, что пара $(u_*,v_*)$ является седловой точкой. Для завершения доказательства заметим, что приведенные рассуждения обратимы.
\epr

Из свойств функции $\varphi(\cdot)$ и пространств $\UB$, $\VB$, требуемых в условии теоремы \ref{teoFC-TfKK}, следует замкнутость графика отображения $\Fg{F}{\varphi}:\UB\times\VB\mapsto2^{\UB\times\VB}$, компактность и хаусдорфовость пространства $\UB\times\VB$ с топологией произведения топологических пространств $\UB$, $\VB$. Таким образом, выполнены условия теоремы \ref{TfO} на отображение $\Fg{F}{\varphi}$. 

Пусть для отображения $\Fg{\barf{F}}{\varphi}$ выполняется условие \fref{con-cov-not-empt}. Тогда из теоремы \ref{TfO} и леммы \ref{lem_fix_eq_saddle}  следует, что множество седловых точек \fref{sedlo} для функции $\varphi$ не пусто. Значит, имеет место равновесие \fref{s1}. 

Обратно, пусть имеет место равновесие \fref{s1}. Тогда из условий на функцию $\varphi(\cdot)$ и пространства $\UB$, $\VB$ следует, что множество
$$
\argmin_{u\in\UB}\max_{v\in\VB}\varphi(u,v)\times\argmax_{v\in\VB}\min_{u\in\UB}\varphi(u,v)
$$
не пусто и совпадает со множеством седловых точек отображения $\Fg{F}{\varphi}$. 
В этом случае из леммы \ref{lem_fix_eq_saddle} и  теоремы \ref{TfO} следует выполнение условия \fref{con-cov-not-empt}.

\subsection{Схема доказательства теоремы \ref{NEQ-TfKK}}

Приведенная ниже лемма следует из определения отображения $F$.

\begin{lemRU}\label{lem_fix_eq_NE}
\beq\label{nesh-eq-fix}
N(J)=\fix{\Fg{F}{(X,J)}}.
\eeq
\end{lemRU}

Из свойств функций $J_i$ и пространств $X_i$, $i\in\nint{1}{n}$, требуемых в условии теоремы \ref{NEQ-TfKK}, следует замкнутость графика отображения $\Fg{F}{(X,J)}:X\mapsto2^{X}$, компактность и хаусдорфовость пространства $X$. Таким образом, для пространства $X$ и отображения $\Fg{F}{(X,J)}$ выполнены условия теоремы \ref{TfO}. 

После этого замечания нетрудно видеть, что первая часть теорема \ref{NEQ-TfKK} следует из леммы \ref{lem_fix_eq_NE} и теоремы  \ref{TfO}. 
Эквивалентность условий \fref{kur-in-O} и \fref{con-cov-not-empt} следует из определения функции $\Fg{\barf{F}}{(X,J)}$ по отображению $\Fg{F}{(X,J)}$.

\end{document}